\documentclass[leqno,11pt]{article}
\usepackage{amssymb,amsmath}
\usepackage{psfig}

\begin{document}

\newcommand{\qed}{\hfill$\blacktriangleright$\medskip}
\newcommand{\proof}[1]{\noindent{\bf Proof#1:\  }}
\newcommand{\mathto}{\mathop{\longrightarrow}\limits}
\def\span{\mathrm {span}}
\def\const{\mathrm {const}}
\def\btl{$\blacktriangleleft$}
\def\bt{\begin{theorem}}
\def\et{\end{theorem}}
\def\bp{\begin{pusto}}
\def\ep{\end{pusto}}
\def\bi{\begin{itemize}}
\def\ei{\end{itemize}}
\def\p{\partial}
\def\e{\varepsilon}
\def\d{\delta}
\def\X{X^{(r)}}
\def\mn{\medskip\noindent}
\def\n{\noindent}
\def\Op{{\cal O}{\it p}\,}
\def\Cal{\cal}
\def\ol{\overline}
\def\A{{\Cal A}}
\def\E{{\Cal E}}
\def\F{{\Cal F}}
\def\G{{\Cal G}}
\def\H{{\Cal H}}
\def\R{{\Cal R}}
\def\P{{\Cal P}}
\def\T{{\Cal T}}
\def\L{{\Cal L}}
\def\S{{\Cal S}}
\def\U{{\Cal U}}
\def\D{{\Cal D}}
\def\V{{\Cal V}}
\def\W{{\Cal W}}
\def\Z{{\Cal Z}}
\def\C{{\Cal C}}
\def\O{{\Cal O}}
\def\BH{{\bold H}}
\def\BF{{\bold F}}
\def\BS{{\bold S}}
\def\bbZ{\mathbb {Z}}
\def\bbR{{\mathbb R}}
\def\bbS{{\mathbb S}}
\def\bbC{{\mathbb C}}
\def\Int{ \mathrm {Int}}
\def\Vert{ \mathrm {Vert}}
\def\Ker{\mathrm {Ker}}
\def\id{\mathrm {id}}
\def\su{\mathrm {supp}\,}
\def\Id{\mathrm {Id}}
\def\top{\mathrm{top}}
\def\bot{\mathrm{bot}}
\def\wt{\widetilde}
\def\inv{\mathrm{inv}}
\def\Diff{\mathrm{Diff}\,}
\def\Hol{\mathrm {Hol}}
\def\Sec{\mathrm{Sec}}
\def\Sol{\mathrm{Sol}}
\newtheorem{theorem}{Theorem}[subsection]
\newtheorem{corollary}[theorem]{Corollary}
\newtheorem{lemma}[theorem]{Lemma}
\newtheorem{proposition}[theorem]{Proposition}
\newtheorem{definition}[theorem]{Definition}
  \newtheorem{remark}[theorem]{Remark}
\newtheorem{example}[theorem]{Example}
\newtheorem{exercise}[theorem]{Exercise}
\newtheorem{conjecture}[theorem]{Conjecture}
\newtheorem{question}[theorem]{Question}
\newtheorem{problem}{Problem}
\newtheorem{pusto}[theorem]{}
\newtheorem{criterion}[theorem]{Criterion}
\newtheorem{generalization}[theorem]{Generalization}


\title
{ \hskip 2.0truein {\it To Andr\'e Haefliger with admiration
}
\bigskip
\bf
      Holonomic approximation and Gromov's h-principle
}

\author  {{\bf Y. M. Eliashberg}\\
     Stanford University, USA \\
\and
 {\bf N. M. Mishachev}\\
    Lipetsk Technical University, Russia \\
    and Stanford University, USA
}
\date{}
\maketitle

\begin{abstract}
\noindent
In 1969 M. Gromov in his PhD thesis \cite{G1} greatly
generalized Smale-Hirsch-Phillips immersion-submersion theory
(see \cite{Sm},\cite{Hi}, \cite{Ph}) by proving  what is now
called the $h$-principle for invariant open differential relations
over open manifolds. Gromov extracted the original geometric idea
of Smale and put it to work in the maximal possible generality.
Gromov's thesis was brought to the West by A. Phillips and was popularized in his talks. However,
most western mathematicians first learned about Gromov's theory from
 A Haefliger's article \cite{H}. The current paper is devoted
to the same subject as the papers of Gromov and Haefliger. It
seems to us that we further purified Smale-Gromov's  original idea by
extracting from it a simple but very general theorem about holonomic
approximation of sections of jet-bundles
(see Theorem \ref{thm:HAT} below). We show below
that Gromov's theorem
as well as some other results in the   $h$-principle spirit are
immediate corollaries of Theorem \ref{thm:HAT}.
\footnote{The authors
are partially supported by the National Science Foundation}
\end{abstract}

\section {Holonomic approximation}\label{s.fph}

\subsection {Jets and holonomy}\label{ss.sjh}

\mn
Given a $C^\infty$-smooth fibration $p:X\to V$, we denote by $\X$ the
space of $r$-jets of smooth sections $f:V\to X$ and by $J^r_f:V\to\X$ the
$r$-jet of a section $f:V\to X$. When the fibration $X=V\times W\to V$ is trivial
then the space $\X$ is sometimes denoted by
$J^r(V,W)$, and called the space of $r$-jets
of maps $V\to W$.
A section $F:V\to \X$ is called {\it holonomic} if it has the form
$J^r_f$ for a section $f:V\to X$. The correspondence $f\mapsto
J^r_f$ defines the derivation map $J^r:\Sec\,X\to \Sec\,\X$. Its
one-to-one image $J^r(\Sec\,X)$ coincides with the space $\Hol\X
\subset \Sec\X$ of holonomic sections, i.e. we have
$$\Sec\,X\,\mathop{\simeq}\limits^{J^r} \,\Hol\,\X\,\hookrightarrow
\,\Sec\,\X\,.$$ Notice that the $C^0$-topology on
$\Sec\,\X$ induces via $J^r$ the $C^r$-topology on $\Sec\,X$.

\mn
Following Gromov's book \cite{G2}   we will denote by $\Op A$
an arbitrary small but non-specified (open) neighborhood of
a subset $A\subset V$.  We will assume that the manifold $V$
and the bundle $\X$ are endowed with   Riemannian metrics and denote
by $U_\e(A)$ the  metric $\e$-neighborhood $\{x|\;{\mathrm {dist}}(x,A)<\e\}$.

\mn
Given an arbitrary subset $A\subset V$  a section $F:A\to\X$
is called {\it holonomic} if there exists a holonomic section
$\widetilde F :\Op A\to \X$ such that $\widetilde F|_A=F$.
A section $F:V\to X$ is called {\it holonomic over} $A\subset V$
if the restriction $F|_A$ is holonomic.    Given a fibration $\pi:V\to B$ we
say that a section $F:V\to\X$ is {\it fiberwise holonomic}
if there exists a continuous family
of holonomic sections $\wt F_b:\Op\pi^{-1}(b)\to \X,\; b\in B$,
such that for each $b\in B$ the sections
$\wt F_b$ and $F$ coincide over the fiber $\pi^{-1}(b)$.
Note the following trivial but important fact:
\bp
\label{lm:point}
Any section
$F:V\to X$ is  holonomic over any point  $v\in V$. Moreover,
it is fiberwise holonomic with respect to the trivial
fibration $\,\Id:V\to V$.
\ep
Indeed,  we can
take the Taylor polynomial map  which corresponds to $F(v)$ with
respect to some local coordinate system centered at $v$
as a section $\wt F_v: \Op v\to \X$. Moreover, the local coordinate system
can be chosen smoothly depending on its center.
\qed
\subsection{Holonomic approximation}

{\bf Question:} {\it Is it possible to approximate any
section $F:V\to\X$ by a holonomic section? In other words, given a
$r$-jet section and an arbitrary small neighborhood of the image of this
section in the jet space, can one find a holonomic section in
this neighborhood?}

\mn
The answer is  evidently negative
(excluding, of course, the situation when the initial
section is already holonomic). For instance, in the case $r=q=1$
the question has the following geometrical reformulation: given a
function and a $n$-plane field along the graph of this function, can
one $C^0$-perturb this graph to make it almost tangent to
the given field?

\mn
The problem of finding a holonomic approximation of a  section
of the $r$-jet  space {\it near a  submanifold}
$A\subset\bbR^n$ is also usually unsolvable. The only exception is
the
zero-dimensional case: as we already stated above in \ref{lm:point}
 any section can be approximated near any point
by the $r$-jet of the respective Taylor polynomial map.

\mn
In contrast, the following theorem says that we {\it always
can  find} a holonomic approximation of a section $F:V\to \X$
{\it near a slightly deformed} submanifold $\wt A\subset V$
if the original  set $A\subset V$ is of positive codimension.

\bp
\label{thm:HAT}
{\bf (Holonomic Approximation Theorem)}
Let $A\subset V$ be a polyhedron of positive codimension
and $F:\Op A\to \X$ be a section. Then for arbitrary small $\;\delta,\e>0$
there exist a
diffeomorphism $h:V\to V$  with $\,||h-\Id||_{C^0}<\delta$
and a holonomic section
$\wt F:\Op h(A)\to \X$ such that
the image
$h(A)$ is contained in the domain of the definition of the section $F$ and
$$||\wt  F- F|_{\Op h(A)}||_{\,C^0}<\e\,.$$
\ep

\mn
 We use the term {\it polyhedron} here in the sense that $A$
is a subcomplex of a certain smooth triangulation of  the
manifold $V$.

\mn
As we will see below, the {\it relative} and the
{\it parametric} versions of the theorem are also true.
In the relative version the section $F$ is assumed to be already  holonomic
over $\Op B$, where $B$ is a subpolyhedron of $A$, while the diffeomorphism
$h$ is constructed to be fixed on $\Op B$ and $\wt F$
is required to coincide with $F$ on $\Op B$.
Here is the parametric version of \ref{thm:HAT}.

\bp
\label{thm:PHAT}
{\bf (Parametric Holonomic Approximation Theorem)}
Let $A\subset V$ be a polyhedron of positive codimension
and $F_z:\Op A\to \X$ be a family of sections
parametrized by a cube $I^m$, $m=0,1,\dots\;$.
Suppose that the sections $F_z$ are holonomic for $z\in\Op \partial I^m$.
Then for arbitrary small $\;\delta,\e>0$
there exist a family of
diffeomorphisms $h_z:V\to V$ and a family of holonomic section
$\wt F_z:\Op h(A)\to \X$, $z\in I^m$,
such that
\begin{itemize}
\item $||h_z-\Id||_{C^0}<\delta$;
\item $h_z=\Id$ for $z\in \Op\partial I^m$;
\item $ ||\wt F_z - F_z|_{\Op h_z(A)}||_{\,C^0}<\e$ ;
\item $\wt F_z=F_z$ for $z\in\Op\partial I^m$.
\end{itemize}
\ep

\mn
Using the induction over skeleta of   the polyhedron $A$, and
taking into account that the fibration $X\to V$ is trivial over simplices,
we reduce Theorem \ref{thm:HAT} to its special case for the pair
$(A,B)=(I^k, \partial I^k)\subset \bbR^n\,.$
We consider this special case in the next section.

\subsection {Holonomic Approximation over a cube}

\bp
\label{thm:FIT}
{\bf (Holonomic Approximation over a cube)}
Let $I^k\subset\bbR^{n}$, $k<n$, be the unit
cube in the coordinate subspace $\bbR^k\subset\bbR^{n}$ of the first
$k$ coordinates. Then for any section
$$F:\Op I^{k}\to J^r(\bbR^{n},\bbR^q)\,$$
which is holonomic over  $\Op \partial I^k$
and for arbitrary small
$\;\delta,\e >0$ there exist a
diffeomorphism $h:\bbR^{n}\to\bbR^{n}$ of the form
$$h(x_1,\dots,
x_{n})=(x_1,\dots,x_{n-1},x_{n}+\varphi(x_1,\dots,x_{k})),$$
and a holonomic section
$$\wt F:\Op h(I^k)\to J^r(\bbR^{n},\bbR^q)$$
such that
\begin{itemize}
\item $|\varphi|<\delta\,;\;\; \varphi|_{\Op\partial I^k}=0\,$;
\item  the image $h(I^k)$ is contained in the domain of the
definition of the section $F$;
\item  $\wt F|_{\Op \partial I^k}=F|_{\Op \partial I^k}\,$ and
\item  $||\wt F-F|_{\Op h(I^k)}||_{C^0}<\e\,.$
\end{itemize}
\ep

\mn
Theorem \ref{thm:FIT} will be deduced from   Inductional
lemma \ref{lm:ind} which we formulate below.

\mn
Let  $\pi_l: I^k\to I^{k-l}$, $l=1,\dots,k,$ be the projection
$$(x_1,\dots,x_k)\mapsto (x_1,\dots, x_{k-l})$$
whose fibers are $l$-dimensional cubes
$$I^l(y)=y\times I^l\,,\, y=(x_1,\dots, x_{k-l})\in I^{k-l}\,.$$
Note that for $l<k$ we have
$I^l(y)=\bigcup\limits_{t\in I}I^{l-1}(y,t)\,$
where $t=x_{k-l+1}$.

\bp
\label{lm:ind}
{\bf (Inductional lemma)}
{\it Suppose that  a section
$$F: \Op I^k\to J^r(\bbR^{n},\bbR^q)$$
is holonomic over $\Op\partial I^k$ and for a positive integer $l\leq k$
fiberwise holonomic with respect to the fibration
$\pi_{l-1}:I^k\to I^{k-l+1}$.
Then for arbitrary small $\;\delta,\e>0 $ there exist
a diffeomorphism $h:\bbR^{n}\to\bbR^{n}$ of the form
$$h(x_1,\dots, x_{n})=(x_1,\dots,x_{n-1},x_{n}+\varphi(x_1,\dots,x_{k}))\,,$$
and a section
$$\wt F:\Op h(I^n)\to J^r(\bbR^{n},\bbR^q)\,,$$
such that

\bi
\item $|\varphi|<\delta\,$; $\;\varphi|_{\Op\partial I^k}=0\,$;
\item the image $h(I^k)$ is contained in the domain of the
definition of the section $F$;
\item $\wt F|_{\Op \partial I^k}=F|_{\Op \partial I^k}\,$;
\item $||\wt F-F|_{\Op h(I^k)}||_{C^0}<\e\,$;
\item the section $\wt F|_{h(I^k)}$ is fiberwise holonomic with respect to
the fibration
$$\pi_l\circ h^{-1}:h(I^k)\to I^{k-l}\,.$$
\ei}
\ep
\mn
\proof{ of Inductional lemma}
We recommend to the reader to keep in mind the simplest case
$n=2\,,\,k=l=1$ while reading the proof for the first time.

\mn
For $(y,t)\in I^{k-l}\times I$ set
$$U_\d(y,t)=U_\d (I^{l-1}(y,t))\;\;{\rm and}\;\;
U^\partial_\d(y,t)=U_\d (\partial I^{l-1}(y,t))\,.$$
It follows from  the definition of a fiberwise holonomic section  and compactness
arguments that we can choose $\d>0$ so small that
there would exist a continuous family of holonomic sections
$$F_{y,t}=J^r_{f_{y,t}}:U_\d(y,t)
\to J^r(\bbR^{n},\bbR^q),
\;\;(y,t)\in I^{k-l+1}=I^{k-l}\times I,$$
such that for all $(y,t)\in  I^{k-l}\times I$:
\bi
\item $F_{y,t}$ is defined on $U_\d(y,t)\,$;
\item $F_{y,t}|_{I^{l-1}(y,t)}=F|_{I^{l-1}(y,t)}$\,.
\ei
We can further adjust $\d$ and the family $F_{y,t}$ in order to have
\bi
\item $F_{y,t}|_{U^\partial_\d(y,t)}
=F|_{U^\partial _\d (y,t)}$ for all
$(y,t)\in  I^{k-l}\times I\,$ and
\item $F_{y,t}=F|_{U_\d (y,t)}$ for all
$(y,t)\in \Op \partial(I^{k-l}\times I)$.
\ei

\mn
For a sufficiently large $N$, which is determined by the
Interpolation Property below, we set
\begin{equation*}
\begin{split}
&U_i(y)= U_\d (y,\frac{i}{2N})\,,\;
U^\partial_i(y)= U^\partial_\d (y,\frac{i}{2N})\,,\; y\in
I^{k-l},\;i=0,\dots,2N\,,\\
&U_i^{\top}(y)=U_i(y)\cap \{x_n\geq\frac{\d}{2}\}\,,\\
&U_i^{\bot}(y)=U_i(y)\cap \{x_n\leq-\frac{\d}{2}\}\,.\\
\end{split}
\end{equation*}
We will write $F^i_{y }$ and $f^i_{y }\,$ instead of
$F_{y,\frac{i}{N}}$ and $f_{y,\frac{i}{N} }$.

\mn
Notice that   $$ \sigma(N) =\max\limits_{ y\in
I^{k-l},\; i=1,\dots, 2N,\;  x\in U_i(y)\cap U_{i-1}(y)}
||F^i_{y}(x)-F^{i-1}_{y}(x)||
\mathop{\to}\limits_{N\to\infty} 0,$$
and hence we  have the following

\mn
{\bf Interpolation Property.}
{\it For any $\e>0$,  a sufficiently large $N$
and all odd integers $i=1,3,\dots 2N-1$ there exist continuous
families of holonomic sections
$$G^i_{y}=J^r_{g^i_{y }}:U_i(y)\to  J^r(\bbR^{n},\bbR^q)\,,\;
y\in I^{k-l},$$
such that
\bi
\item$G^i_y$ interpolate between
$F^{i-1}_y$ and $\,F^{i+1}_y\,$ for all $y\in I^{k-l}\,:$
$$G^i_{ y}=
\begin{cases}
& F^{i+1}_{y} \;\;\hbox{on}\;\;  U_{i}^{\top} (y)\cap U_{i+1}^{\top}\,,\\
& F^{i-1}_{y}\;\;\hbox{on}\;\;  U_{i}^{\bot} (y)\cap U_{i-1}^{\bot}\,;\\
\end{cases}$$
\item $||G^i_{y }-F^i_{y }||_{C^0}<\e\,$ for all $y\in I^{k-l}\,;$
\item $G^i_{y}|_{U^\partial_i(y)}
=F^i_y|_{U^\partial_i(y)}=F|_{U^\partial_i(y)}\,$ for all $y\in I^{k-l}\,;$
\item $G^i_{y}|_{U_i(y)}=F^i_y|_{U_i(y)}=F|_{U_i(y)}\,$
for all $y\in \Op \partial I^{k-l}\,.$
\ei}

\bigskip
\mn
For even values $i=0,2,\dots,2N$ we set $G^i_y\equiv F^i_y$.
Take a cut-off function $\theta_N:\bbR\to I$, which is equal to
$0$ on $\Op(\bbR\setminus I)$ and is equal to $1$ on
$[\frac{1}{4N},1-\frac{1}{4N}]$, and define a function
$\varphi_N:\bbR^k\to\bbR$ by the formula
\begin{equation*}
\varphi_N(y,t,x)= \theta_N(t)\theta_N(||x||) \theta_N(||y||)
\cos 2N\pi t,\;\; y\in \bbR^{k-l}\,, t\in \bbR\,, x\in \bbR^{l-1}\,.
\end{equation*}
Consider the map $h:\bbR^n\to\bbR^n$,
$$h(x_1,\dots ,x_{n})=(x_1,\dots,x_{n-1},x_{n}+
\d_1\varphi_N(x_1,\dots,x_{k}))$$
where $\d/2<\d_1<\d$.
Viewing the image $\tilde I^l(y)=h(I^l(y))$ as the union
$$\tilde I^l(y)=\bigcup\limits_0^{2N-1}\wt I_i(y)\,,$$
where
$$\tilde I_i(y)=\tilde I^l(y)\cap\{\frac{i}{2N}\leq t=x_{k-l+1}\leq \frac{i+1}{2N}\}$$
we define a continuous family  of holonomic sections
$$\wt F_y: \Op \tilde I^l(y) \to J^r(\bbR^n,\bbR^q)\,,\;y\in I^{k-l}$$
by the formula
$$\wt F_y|_{\Op \tilde I_i}= G_y^i|_{\Op \tilde I_i} \,,\;\; i=0,\dots,2N-1.$$
Then the section $\wt F:\Op h(I^k)\to J^r(\bbR^n,\bbR^q)$ defined by the formula
$$\wt F(y,t,x)=
\wt F_y(t,x)\;\;\hbox{for}\;\; (y,t,x)\in \Op(I^{k-l}\times I\times I^{n-k-l-1})$$
is holonomic with respect to the fibration
$\pi_l\circ h$ and is $\,\e$-close to $F$.
\qed

\proof{ of Theorem \ref{thm:FIT}}
We will prove the theorem by an induction over $l$.
Consider for $l=0,\dots k\;$ the following

\mn
{\bf Inductional hypothesis $\A^{(l)}$.}
{\it Let $F:\Op I^k\to J^r(\bbR^n,\bbR^q)$
be a section which is holonomic over $\Op\partial I^k$.
For arbitrary small $\d,\e>0$ there exist
a diffeomorphism $h:\bbR^{n}\to\bbR^{n}$
of the form
$$h(x_1,\dots, x_{n})=(x_1,\dots,x_{n-1},x_{n}+\varphi(x_1,\dots,x_{k}))\,,$$
and a section
$$\wt F^l:\Op h(I^n)\to J^r(\bbR^{n},\bbR^q)\,$$
as required by the Inductional lemma, i.e
\bi
\item $|\varphi|<\d\,$; $\;\varphi|_{\Op\partial I^k}=0\,$;
\item the image $h(I^k)$ is contained in the domain of the
definition of the section $F$;
\item $\wt F^l|_{\Op \partial I^k}=F|_{\Op \partial I^k}$;
\item $||\wt F^l-F|_{\Op h(I^k)}||_{C^0}<\e$;
\item the section $\wt F^l|_{h(I^k)}$ is fiberwise holonomic with respect to
the fibration
$$\pi_l\circ h^{-1}:h(I^k)\to I^{k-l}\,.$$
\ei}

\mn
According to \ref{lm:point} the given  section
$F:I^n\to J^r(\bbR^{n},\bbR^q)$,
which is holonomic near $\partial I^k\,$, is tautologically fiberwise
holonomic with respect to the fibration  by points
$\pi_0:I^k\to I^k$.
This implies $\A^{(0)}$ and thus gives us the base for the induction.
For $l=1$ the implication  $\A^{(l-1)}\Rightarrow \A^{(l)}$ follows immediately
from Inductional lemma \ref{lm:ind}, but in the general case $l>1$
we cannot apply \ref{lm:ind} directly    because the section $F^{l-1}$
is defined near the deformed cube rather than the original one.
Notice  however, that the diffeomorphism $h:\bbR^n\to\bbR^n$
induces the covering map
$h_*:J^r(\bbR^{n},\bbR^q)\to J^r(\bbR^{n},\bbR^q)$. The section
$\bar F^{l-1}=(h_*)^{-1}(\wt F^{l-1})$ is defined over
$\Op I^k\,$, coincides with $F$ near $\partial I^k$ and
fiberwise holonomic with respect to
the fibration $\pi_{l-1}:I^k\to I^{k-l+1}$.
Applying  Inductional lemma \ref{lm:ind}
we can approximate $\bar F^{l-1}$ by
a section $\wt F'$ over a deformed cube $h'(I^k)$. The  section $\wt F'$  coincides
with $\bar F^{l-1}$ near $\partial I^k$ and fiberwise holonomic with respect to
the fibration $\pi_{l}\circ h':h'(I^k)\to I^{k-l}$.
If $\wt F'$ is sufficiently $C^0$-close to $\bar F^{l-1}$, then
the section $\wt F^{l}=h_*(\wt F')$ is the required approximation
of $F$ in a neighborhood of $h''(I^k)$, where $h''=h\circ (h')$.
This proves $\A^{(l)}$ and Theorem \ref{thm:FIT}.
\qed

\paragraph{Parametric case.}
It turns out that  Inductional lemma \ref{lm:ind} implies
also the   parametric version of  Theorem \ref{thm:FIT}.
Namely, we have

\bp
\label{thm:PFIT} {\bf (Parametric version of Theorem
\ref{thm:FIT})}
Let $F_z,\;z\in I^m$, be a family of sections $\,\Op
I^k\to J^r( \bbR^n,\bbR^q)$, parameterized by the cube $I^m$.
Suppose that $k<n$ and the sections $F_z$ are holonomic over
$\Op\partial I^k$ for all $z\in I^m$, and holonomic over the whole
$I^k$ for $z\in\Op \partial  I^m$.
Then for arbitrary small
$\;\delta,\e >0$ there exist a family
of diffeomorphisms $h_z:\bbR^{n}\to\bbR^{n}$ of the form
$$h_z(x_1,\dots,
x_{n})=(x_1,\dots,x_{n-1},x_{n}+\varphi_z(x_1,\dots,x_{k})),$$
and a family of holonomic sections
$$\wt F_z:\Op h_z(I^k)\to J^r(\bbR^{n},\bbR^q)$$
such that
\begin{itemize}
\item $|\varphi_z|<\delta;\;\;\; \varphi_z|_{\Op\partial
I^k}=0; $
\item the section  $F_z$ is defined  in a neighborhood of $h_z(I^k)\; z\in I^m$;
\item  $\wt F_z|_{\Op \partial I^k}=F_z|_{\Op \partial I^k}$
\item $\wt F_z =F_z $ for $z\in\Op I^m$ and
\item  $||\wt F_z-F_z|_{\Op h_z(I^k)}||_{C^0}<\e\,.$
\end{itemize}
\ep

\proof{}
Consider the cube $I^{k+m}=I^k\times I^m\subset
\bbR^n\times\bbR^k=\bbR^{n+m}$. Let $  J^r(\bbR^{n+m}|\bbR^n,\bbR^q)$
be the bundle whose restriction to $\bbR^n\times z,\, z\in\bbR^m,$ equals
$J^r(\bbR^n,\bbR^q)$. The family of sections
$$F_z:I^k\to J^r(\bbR^n,\bbR^q)$$
can be viewed as a section
$$\overline F:I^{k+m}\to J^r(\bbR^{n+m}|\bbR^n,\bbR^q)\,.$$
The section $\overline F$ lifts
to a section $\overline{\overline F}:I^{k+m}\to
J^r(\bbR^{n+m},\bbR^q)$, so that $\pi\circ \overline{\overline
F}=\overline F$, where
$$\pi:J^r(\bbR^{n+m},\bbR^q)\to J^r(\bbR^{n+m}|\bbR^n,\bbR^q)$$
is the canonical projection.
Moreover, the section $\overline{\overline F}$ can be
chosen holonomic near $\partial I^{k+m}$. Hence we can apply
Theorem \ref{thm:FIT} to get an  $\e$-approximation
$\wt{\wt F}$ of $\overline{\overline F}$ over a
$\d$-displaced  cube $h(I^{k+m})$. Then the composition
$\wt F=\pi\circ\wt{\wt F}:I^{k+m}\to
J^r(\bbR^{n+m}|\bbR^n,\bbR^q)$ can be interpreted as the required
family $\{\wt F_z\}_{z\in I^m}$ of holonomic  $\e$-approximations
of the family  $\{F_z\}$ near $\{h_z(I^k)\}$.
\qed

\mn
In the same way  as Theorem \ref{thm:FIT} implies
Theorem \ref{thm:HAT}, i.e. via the induction over skeleta,
Theorem \ref{thm:PFIT} implies Theorem \ref{thm:PHAT}.

\section{Applications}
\subsection{Gromov's h-principle for $\Diff V$-invariant differential
relations over open manifolds}

\paragraph{Differential relations.}
A differential relation (or condition) imposed on sections $\varphi:V\to X$
of a fiber bundle $X\to V$
is a subset $\R\subset \X $, where $r$ is called the {\it order} of $\,\R$.
A  section $\Phi:V\to\X$ is called a
{\it formal solution} of $\R$ if $\Phi(V)\subset\R$.
A $C^r$-section $\varphi:V\to X$ is called {\it a solution} of the
relation $\R$  if $J^r_\varphi(V)\subset \R$.
We will denote the space of solution of $\R$
by $\Sol\,\R$, the space of formal solutions by $\Sec\,\R$,
and the space of holonomic formal solutions by $\Hol\,\R$.
The $r$-jet extension establishes  a one-to-one
correspondence $$J^r:\Sol\,\R\to \Hol\,\R$$
between the solutions and holonomic formal solutions
of $\R$. We will use the term ``solution" also for  the sections from $\Hol\,\R$
when the distinction between the solutions of $\R$ as sections of
$X$ or $\X$ is clear from the context, or irrelevant.

\mn
\paragraph{$\Diff\,V\,$-invariant differential relations.}
Given a fibration $p:X\to V$ we will denote by $\Diff_V X$ the
group of fiber-preserving diffeomorphisms $h_{X}:X\to X$, i.e.
$h_X\in\Diff_VX$ if and only if there exists a diffeomorphism
$h_V:V\to V$ such that $p\circ h_X=h_V\circ p$. Let
$\pi:\Diff_V X\to \Diff V$ be the projection $h_X\mapsto h_V$.
We are interested when this arrow can be reversed, i.e. {\it
when there exists a homomorphism $j:\Diff V\to \Diff_V X$ such
that $\pi\circ j=\id\,$.\,} We call a fibration $X\to V$, together
with a homomorphism $j$, {\it natural} if such a lift exists.
For instance, the trivial fibration $X=V\times W\to V\to  V$ is natural. Here
$j(h_V)=h_V\times\id$. The tangent bundle $T(V)\to V$ is natural as well.
If a fibration $X\to V$ is natural then any fibration
associated with it is natural as well.
In particular, if $X\to V$ is natural then $\X\to V$ is
natural.  The implied lift
$$j^r:\Diff V\to \Diff_V \X\,,\;\; h\mapsto h_*$$
is defined here by the formula
$$h_*(s)=J^r_{j(h)\circ \bar s}(h(v))$$
where $s\in \X\,$, $v=p^r(s)\in V\,$, and $\bar s\,$ is a
local section near $v$ which represents the $r$-jet $s$.

\mn
Given a  natural fibration $X\to V\,$, a differential relation
$\R\subset \X$ is called $\Diff V$-{\it invariant} if  the action
$s\mapsto h_*s,\, h\in\Diff V$ leaves $\R$ invariant.
In other words, a differential relation $\R$  is $\Diff V$-invariant
if it can be defined in $V$-coordinate free form. Notice that
though the definition of a $\Diff V$-invariant relation depends on the
choice of the  homomorphism $j$ this choice is fairly obvious
in most interesting examples, and we will not specify it.

\mn
The action $s\mapsto h_*s$ preserves the set of holonomic sections:
$$h_*(J^r_f)=J^r\left(j(h)\circ f\circ h^{-1}\right),$$
$f\in \Sec X$, $h\in \Diff V$. In particular,
the group $\Diff V$ acts on the space $\Sol\,\R$ of an
invariant differential relation $\R$.

\paragraph{Homotopy principle.}
Existence of a formal solution is a necessary condition
for the solvability of a differential relation $\R$, and thus before
trying to solve $\R$ one should check whether  $\R $ admits a formal
solution. The problem of finding formal solutions is of pure
homotopy-theoretical nature. This problem maybe simple, or  highly non-trivial, but
in any case it is important to treat the homotopical problem
first, and  look for   genuine solutions only  after existence of formal solutions
have been  already  established.

\mn
Though finding a formal solution is  an algebraic, or
homotopy-theoretical problem which is a dramatic simplification
of the original differential problem, and at first thought the existence of a formal solution
cannot be sufficient for the genuine solvability of $\R$, in the second half of XX century
there were discovered    many
large and geometrically interesting classes of differential relations $\R$
for which the solvability of  the formal problem turned out to be  sufficient
for the genuine solvability.  Moreover, in most of these examples
the spaces  of   formal and genuine solutions
are much closer related that one could expect.

\mn A differential relation $\R$  is said to satisfy
the  {\it (parametric) $h$-principle},
or {\it homotopy principle} if the inclusion
$\Hol\,\R\hookrightarrow \Sec\,\R$ is a weak homotopy equivalence.
This means, in particular, that  any formal solution can be
deformed inside $\F$ into a genuine  solution and any two
solutions  which can be joined by a family of formal solutions
can be also joined  by a family of genuine solutions. Moreover,
similar properties hold  for families of formal
and genuine solutions, depending on arbitrary number of parameters.

\mn
In  fact, it is useful to consider   different ``degrees"
of $h$-principles, when one want to establish closer and closer
connection between formal and genuine solutions. For instance,
different forms of $h$-principles may include some approximation
and  extension properties.

\paragraph{Open manifolds.}
A manifold $V$ is called {\it open}, if  there are no
closed manifolds among its connected components. As it is well known
\bp
\label{lm:compression}
If $V$ is open, then there exists a polyhedron $K\subset V$,
$\,codim\,V_0\geq 1$, such that $V$ can be compressed by an
isotopy $g^t:V\to V,\,\,t\in [0,1],$ into an arbitrary small
neighborhood $U$ of $\,V_0$. \ep

\mn
Here is the main theorem of this section.
\bp
\label{thm:open}
{\bf (Gromov, 1968)} Let $V$ be an open manifold and
$X\to V$ a natural fiber bundle.
Then any open $\Diff V$-invariant differential relation $\R\subset\X$
satisfies the parametric $h$-principle.
\ep

\proof{}
We need to proof that given a family of sections
$F_z\in\Sec\,\R, \;z\in I^m$, $m=0,\dots$, such that for
$z\in\Op\partial I^m$ the section $F_z$ is holonomic,
there exists a family $\wt F_z\in\Hol\,\R$
which is homotopic in $\Sec\,\R$ to the family $F_z$, $z\in I^m$,
relative to $\partial I^m$.
Let $K\subset V$ be a polyhedron of positive
codimension, as in \ref{lm:compression}.
According to Theorem \ref{thm:PHAT} there exist
a $C^0$-small family of diffeomorphisms $h_z:V\to V$
and a family $G_z\in\Hol_{\Op h_z(K)}\R$, such that
$h_z=\Id $ for $z\in \partial I^m$, and $G_z$ is $C^0$-close
to $F_z$  over $\Op h_z(K)$. Let $U$ be a small neighborhood of $K$,
such that for all $z\in I^m$ the image $h_z(U)$ is contained in
the domain of the definition of the section $G_z$.
Let $g^t:V\to V,\; t\in [0,1]$, be the isotopy
compressing $V$ into the neighborhood $U$. Then for every $z\in I^m$
the isotopy $g^t_z=h_z\circ g^t\circ h_z^{-1}\;, t\in[0,1],$
compresses $V$ into a neighborhood of $h_z(K)$ where the section
$G_z$ is defined. The desired family of sections
$\wt F_z\in\Hol\,\R$ can be now defined by the formula
$$\wt F_z=(g^1_z)_*^{-1}(G_z),\;z\in I^m\,,$$
where  $(g^1_z)_*:\X\to\X$ is the induced action of
the diffeomorphism $g^1_z$ on the natural fibration $\X\to V$.
\qed

\n
Notice, that   the   {\it $C^0$-dense} $h$-principle
for $\R$ does not necessarily hold, though
the $C^0$-dense $h$-principle holds for
$\R$ {\it near} $A$. The relative $h$-principle does not hold
either. However, we have its following   version.

\bp
Let $\R\subset \X$ be an open $\Diff V$-differential relation
over an open manifold $V$. Let $B\subset V$ be a closed subset
such that each connected component of the complement $V\setminus B$
has  a boundary point which is not in $B$. Then the relative parametric
$h$-principle holds for $\R$ and the pair $(V,B)$.
\ep

\subsection{Directed embeddings of open manifolds}\label{ss.dir}

Let $Gr_n(W)$ be the Grassmanian bundle of tangent $n$-planes
to a $q$-dimensional manifold $W$, $q>n$, and $A\subset Gr_n(W)$ be an arbitrary subset.
 A map $f:V\to W$ of a $n$-dimensional manifold $V$ is called {\it A-directed}
if the tangential (Gauss) map $G_f:V\to Gr_n(W)$ sends $V$ into
$A$. Any $A$-directed map is automatically an immersion. If $A$ is
{\it open}, then the $h$-principle for $A$-directed maps $V\to W$
follows immediately from \ref{thm:open}, since the corresponding
differential relation $\R_A$ is open and $\Diff\,V$-invariant. For
$A$-directed {\it embeddings} Gromov  proved in \cite{G2} via his convex
integration technique  the following \footnote{
Gromov's proof is discussed in \cite{Sp}.
Rourke and Sanderson gave two independent proofs of this theorem in
\cite{RS1} and \cite{RS2}. The proof in \cite{RS2} is quite close to ours.}

\bp
\label{thm:dir}
{\bf (Homotopy principle over embeddings)}
If $\;A\subset Gr_n(W)\;$ is open and $V$ is an open manifold, then
every embedding $f_0:V\to W$, whose tangential lift
$\;G_0=G_{f_0}:V\to Gr_n(W)$ is homotopic over embeddings to a map
$G_1:V\to A\subset Gr_n(W)$ can be
isotoped to a  $A$-directed embedding $f_1:V\to W$.
\ep

\mn
Here the homotopy {\it over embeddings} means that the
underlying homotopy $g_t:V\to W$ is an isotopy.
We show below how Theorem \ref{thm:dir} can be deduced
from Holonomic approximation theorem  \ref{thm:HAT}.

\mn
\proof{}
Suppose that the homotopy
$G_t$ is fixed on the base (that is the underlying
isotopy $g_t:V\to W$ is fixed on $V$) and {\it small} in the following
sense: the angle between $T_v(V)$ and
$G_t(T_v(V))$ is less then, say, $\pi/4$ for all $v\in V$
and $t\in [0,1]$. Let $K\subset V$ be a codimension
$\geq 1$ compact subcomplex in $V$, such that $V$ can be
compressed into an arbitrary small neighborhood of $K$.
The plane fields $G_t(T(V))$ along $V$ define a homotopy of section
$F_t:V\to X^{(1)}$ of the fibration
$X^{(1)}\to V$, where $X$ is a tubular neighborhood of $V\subset W$.
Using  Theorem \ref{thm:HAT}  we can construct
a holonomic approximation $\widetilde F_1$ of
$F_1$ near $h_1(K)$, where
$h_t:V\to V$ is a diffeotopy. For a sufficiently close approximation
the section $\widetilde F_1$ over $h_1(K)$ still lies in (the open set) $A$.
Therefore the graph of the $0$-jet part of this new section over
a neighborhood of $h_1(K)\subset V$ is the image of the required
embedding $f_1:V\to W$. The respective isotopy $f_t$ is the composition of
the compression and the isotopy of sections.

\mn
Any fixed on $V$ (non-small) homotopy $G_t$
may be decomposed into a finite sequence of some small homotopies.
Therefore we can apply the previous construction successively
and get the required isotopy $f_t$.

\mn
If the  general case when the homotopy $G_t$ is not fixed on $V$
we may apply the previous construction to the  pull-back
homotopy $(d\varphi_t)^{-1}\circ G_t$ and the  pull-back
relation $\R_{\widetilde A}= \varphi^*_1(\R_A)$, where
$\varphi_t: W\to W$ is a diffeotopy, which extends the
 underlying for $G_t$  isotopy $g_t:V\to W$.
\qed

\mn
The parametric $h$-principle  is  also valid in this case (with the
same proof).

\subsection{Final remarks}

\mn
Geometrically interesting examples covered by Theorem \ref{thm:open}
include (see \cite{G2}) in the case of open manifolds immersions,
submersions, $k$-mersions (i.e. mappings of rank $\geq k$),
mappings with non-degenerate higher-order osculating spaces,
mappings transversal to foliations, or more generally, to arbitrary
distributions, construction of generating system of exact differential
$k$-forms, symplectic and contact structures
on open manifolds, etc.

\mn On the other hand, the {\it micro-extension} trick which goes
back to M. Hirsch,  allows sometimes to reformulate the  problems
about closed manifolds in terms  of open manifolds.
For example, if $\dim W>\dim V$ then a construction of an
immersion $V\to W$ homotopic to a map $f:V\to W$ is equivalent to
a construction of  an immersion $E\to W$, where $E$ is the total space
of the normal bundle to $T(V)$ in $f^*T(W)$. The manifold $E$ is
open and hence the $h$-principle  \ref{thm:open} applies. The same
trick yields in some cases the directed embedding theorem for
closed manifolds. For instance, Theorem \ref{thm:dir} allows, when
the necessary homotopical conditions are met, to perturb a closed
$n$-dimensional   submanifold $V\subset \bbR^q$  via an isotopy so
that  its projection to $\bbR^k$, where $n<k<q$, becomes an
immersion.

\mn
Gromov  showed in [G1] that the  relative parametric $h$-principle is
equivalent to the relative $h$-principle for foliated manifolds. For
instance, a manifold $V$ with a $k$-dimensional foliation $\F$, such that
the tangent to the foliation  bundle $\xi=T(\F)$ is trivial, admits a
map to $\bbR^{k+1}$ whose restriction to each leaf is an
immersion. A similar statement for a non-integrable sub-bundle
$\xi\subset T(V)$ does not  follow from
Theorem  \ref{thm:open}. To prove this and some other geometric corollaries
one needs  several modifications of the $h$-principle
\ref{thm:open}, which can be found in Gromov's book \cite{G2}.
One of these modifications is concerned with the differential relation
invariant with respect to certain subgroups of $\Diff V$. For instance, the
statement about mappings non-degenerate along
(not necessarily integrable)  tangent distributions
requires the following theorem from Gromov's book \cite{G2}.

\bp
\label{prop:1dim}
Let $X\to V\times\bbR$ be a natural fibration and $\R\subset \X$
be an open  differential relation invariant with respect to
diffeomorphisms of the form
$$(x,t)\mapsto (x,h(t))\;\;x\in V,t\in\bbR.$$
Then $\R$ satisfies the parametric $h$-principle.
\ep

\mn
Proposition \ref{prop:1dim} is also a corollary of our Theorem
\ref{thm:FIT}. The proof follows the same scheme as the proof of
\ref{thm:open} with an additional remark that the perturbation
$h$ implied by \ref{thm:FIT} has a special form precisely as
Proposition \ref{prop:1dim} requires.

\mn
Some other important geometric applications in Gromov's book, e.g.
symplectic and contact embeddings theorems,
need  a generalization of Theorem \ref{thm:open}
from open to   {\it micro-flexible} differential
relations. We will not discuss here the notion of
micro-flexibility and only note that this property is equivalent
to the Interpolation Property considered in the proof of Theorem
\ref{thm:FIT}. Hence it is straightforward to adjust Theorem \ref{thm:HAT}
so that it would serve the needs of the micro-flexible case.
We refer the reader to out forthcoming book \cite{EM} for the
details and further applications of Holonomic Approximation
Theorem.


{\small
\begin{thebibliography}{999}

\bibitem[EM]{EM} Y. Eliashberg and N.M. Mishachev,
{\it Flexible integration. Introduction to h-principle}, in preparation

\bibitem[G1]{G1} M. Gromov, {\it Stable maps of foliations into manifolds,}
Izv. AN SSSR, ser. mat. {\bf 33}(1969), 707--734.

\bibitem[G2]{G2} M. Gromov, {\it Partial differential relations},
Springer-Verlag, 1986


\bibitem[H]{H} A. Haefliger, {\it Lectures on the theorem of Gromov},
Lecture Notes in Math., Vol. {\bf 209}(1971), pp. 128--141

\bibitem[Hi]{Hi} M. Hirsch, {\it Immersions of manifolds},
Trans. Amer. Math. Soc. {\bf 93}(1959), 242--276

\bibitem[Ph]{Ph} A. Phillips, {\it Submersions of open manifolds},
Topology {\bf 6}(1967), 171--206

\bibitem[RS1]{RS1} C. Rourke, B. Sanderson, {\it The compression theorem}
(1997), preprint
\bibitem[RS2]{RS2} C. Rourke, B. Sanderson, {\it Directed
embeddings: a short proof of Gromov's theorem} (2000), preprint

\bibitem[Sm]{Sm} S. Smale, {\it The classification of immersions of spheres
 in Euclidean spaces}, Ann. of Math. (2) {\bf 69}(1959), 327--344

\bibitem[Sp]{Sp} D. Spring, {\it Directed embeddings and the simplification
of singularities} (2000), preprint

\end{thebibliography}}
\end{document}